\documentclass[11pt,draft]{article}%

\title{Multiplier Hopf algebras imbedded in locally compact quantum groups}%
\author{K. De Commer\footnote{Research Assistant of the Research Foundation - Flanders (FWO -
Vlaanderen).}\and A. Van Daele}%
\date{}%

\usepackage{mathrsfs}%
\usepackage{color}%
\usepackage{amssymb, amsfonts, amsthm, amsxtra, amsmath}%
\usepackage{latexsym}%

\newtheorem{Theorem}{Theorem}[section]%
\newtheorem{Lemma}[Theorem]{Lemma}%
\newtheorem{Proposition}[Theorem]{Proposition}%

\newtheorem{Cor}[Theorem]{Corollary}
\newcommand{\Ass}{\textbf{\textit{Assumption: }}}

\begin{document}
\maketitle

\abstract{\noindent Let $(A,\Delta)$ be a locally compact quantum
group and $(A_0,\Delta_0)$ a regular multiplier Hopf algebra. We
show that if $(A_0,\Delta_0)$ can in some sense be imbedded in
$(A,\Delta)$, then $A_0$ will inherit some of the analytic structure
of $A$. Under certain conditions on the imbedding, we will be able
to conclude that $(A_0,\Delta_0)$ is actually an algebraic quantum
group with a full analytic structure. The techniques used to show
this, can be applied to obtain the analytic structure of a
$^*$-algebraic quantum group \textit{in a purely algebraic fashion}.
Moreover, the \textit{reason} that this analytic structure exists at
all, is that the one-parameter groups, such as the modular group and
the scaling group, are diagonizable. In particular, we will show
that necessarily the scaling constant $\mu$ of a $^*$-algebraic
quantum group equals $1$. This solves an open problem posed in
\cite{Kus1}.}

\section*{Introduction}

In \cite{VDae1}, the second author introduced \textit{multiplier
Hopf algebras}, generalizing the notion of a Hopf algebra to the
case where the underlying algebra is not necessarily unital. In
\cite{VDae2}, he considered those multiplier Hopf algebras that
allow a non-zero left invariant functional. It turned out that these
objects, termed \textit{algebraic quantum groups}, possess a very
rich structure, allowing for example a duality theory. These objects
seemed to form an algebraic model of locally compact quantum
groups, which at the time had no proper definition.\\

\noindent In \cite{Kus1}, Kustermans showed that a $^*$-algebraic
quantum group (with a positivity condition on the left invariant
functional) naturally gives rise to a \textit{C$^*$-algebraic
quantum group}, which by then had been defined in a preliminary way
by Woronowicz, Masuda and Nakagami. Kustermans showed however that
there was one discrepancy with the proposed definition, in that the
invariance of the scaling group
with respect to the left Haar weight was only relative.\\

\noindent These investigations culminated in the by now accepted
definition of a \textit{locally compact quantum group} by Kustermans
and Vaes, as laid down in \cite{Kus2}. This definition was (up to
the relative invariance of the scaling group) equivalent with the
one proposed by Woronowicz, Masuda and Nakagami, but the set of
axioms was smaller and simpler. These axioms were very much inspired
by those of $^*$-algebraic quantum groups, but introducing analysis
made it much harder to show that they were sufficiently powerful to
carry a
theory of locally compact quantum groups with the desired properties.\\

\noindent \textit{In this article,} we examine a converse of the
problem studied in \cite{Kus1} and \cite{Kus5}. Namely, instead of
starting with a $^*$-algebraic quantum group and imbedding it into a
locally compact quantum group, we start with an imbedding of a
general regular multiplier Hopf algebra in a locally compact quantum
group, and look whether the multiplier
Hopf algebra inherits some structural properties.\\

\noindent The study of this problem led us to an enhanced structure
theory for $^*$-algebraic quantum groups. For example, the analytic
structure of these objects is a consequence of the fact that all the
actions at hand are diagonizable. This has as a nice corollary that
the scaling constant of a $^*$-algebraic quantum group is
necessarily 1. It is odd that $^*$-algebraic quantum groups, which
provided a motivation for allowing relative invariance of the
scaling group under the Haar weight, turn out to have proper invariance after all.\\

\noindent The paper is organized as follows. \textit{In the first
part}, we introduce the definitions of the objects
at play and introduce notations.\\

\noindent \textit{In the second part} we investigate the following
problem: if a multiplier Hopf algebra $A_0$ can be imbedded in a
locally compact quantum group, does this give us information about
the multiplier Hopf algebra? Firstly, we must specify what we mean
by `imbedded in': $A_0$ has to be a subalgebra of the locally
compact quantum group, and the respective comultiplications
$\Delta_0$ and $\Delta$ have to satisfy formulas of the form
$\Delta_0(a)(1\otimes b)=\Delta(a)(1\otimes b)$ for $a$, $b$ in
$A_0$. Secondly, we must specify whether we imbed $A_0$ in the von
Neumann algebra $M$ or in the C$^*$-algebra $A$ associated to the
locally compact quantum group. Already in the first situation, the
objects of $A_0$ will behave nicely with respect to analyticity of
the various one-parameter groups. But only in the second case can we
conclude, under a mild extra condition, that $A_0$ is invariant
under these one-parameter groups. Moreover, $A_0$ will then
automatically
have the structure of an algebraic quantum group.\\

\noindent\textit{In the third part} we apply the techniques of the
previous section to obtain structural properties of $^*$-algebraic
quantum groups. We want to stress that this section is entirely of
an algebraic nature. For example, we prove in a purely algebraic
fashion the existence of a positive right invariant functional on
the $^*$-algebraic quantum group. Up to now, some involved analysis
was
necessary to arrive at this.\\

\noindent \textit{In the fourth part} we consider some special
cases. We also look at a concrete example. Since examples of
non-compact locally compact groups and $^*$-algebraic quantum groups
are still hard to find, we have to limit ourselves to some already
well-known ones. We will see what can be said about the discrete
quantum group $U_q(su(2))$ in this
context.\\

\noindent Some of the motivation for this paper comes from
$\cite{Lan1}$, where similar questions are investigated in the
commutative and co-commutative case. For example, it is shown that
the function space $C_0(G)$ of a locally compact group contains a
dense multiplier Hopf $^*$-algebra, if and only if $G$ contains a
compact open subgroup. The multiplier Hopf $^*$-algebra will be the
space spanned by translates of regular (=polynomial) functions on
this compact group.

\section{Preliminaries}

In this article, we use the concepts of a regular multiplier Hopf
\linebreak ($^*$-)algebra, a ($^*$-)algebraic quantum group, a
(reduced) C$^*$-algebraic quantum group and a von Neumann-algebraic
quantum group, as introduced respectively in \cite{VDae1},
\cite{VDae2}, \cite{Kus2} and \cite{Kus3} (see also \cite{VDae4}).
Since these objects stem from quite different backgrounds, we will
give a brief overview of their definitions.\\

\textit{Regular multiplier Hopf ($\,^*$-)algebras} \\

\noindent We recall the notion of the \emph{multiplier algebra of an
algebra}. Let $A$ be a non-degenerate algebra (over the field
$\mathbb{C}$), with or without a unit. The non-degeneracy condition
means that if $ab=0$ for all $b\in A$, or $ba=0$ for all $b\in A$,
then $a=0$. As a set, the multiplier algebra $M(A)$ of $A$ consists
of couples $(\lambda,\rho)$, where $\lambda$ and $\rho$ are linear
maps $A\rightarrow A$, obeying the following law:
\[a\lambda(b)=\rho(a)b, \qquad\textrm{for all }a,b\in A.\] In
practice, we write $m$ for $(\lambda,\rho)$, and denote $\lambda(a)$
by $ma$ and $\rho(a)$ by $am$. Then the above law is simply an
associativity condition. With the obvious multiplication by
composition of maps, $M(A)$ becomes an algebra, called the
multiplier algebra of $A$. Moreover, if $A$ is a $^*$-algebra,
$M(A)$ also carries a $^*$-operation: for $m\in M(A)$ and $a\in A$,
we define $m^*$ by $m^*a=(a^*m)^*$ and $am^*=(ma^*)^*$. Note that,
when $A$ is actually a C$^*$-algebra, one may like to impose some
extra continuity conditions on the above linear maps, but actually,
the continuity comes for free, as an application of the closed graph
theorem shows.

There is a natural map $A\rightarrow M(A)$, letting an element $a$
correspond with left and right multiplication by it. Because of
non-degeneracy, this ($^*$-)algebra morphism will be an injection.
In this way, non-degeneracy compensates the possible lack of a unit.
Note that, when $A$ is unital, $M(A)$ is equal to $A$.

Let $B$ be another non-degenerate ($^*$-)algebra. In our paper, a
\emph{morphism} between $A$ and $B$ is a non-degenerate
($^*$-)algebra homomorphism $f:A\rightarrow M(B)$. In general, the
non-degeneracy of a map $f$ means that $f(A)B=B$ and $Bf(A)=B$. If
$A$ and $B$ are both C$^*$-algebras, non-degeneracy means that
$f(A)B$ is norm-dense in $B$. In any case, if $f$ is a morphism from
$A$ to $B$, then $f$ can be extended to a unital ($^*$-)algebra
morphism from $M(A)$ to $M(B)$. A
\emph{proper morphism} between $A$ and $B$ is a morphism $f$ such that $f(A)\subset B$. \\

\noindent We can now state the definition of a regular multiplier
Hopf ($^*$-)algebra. It is the appropriate generalization of a Hopf
($^*$-)algebra to the case where the underlying algebra need not be
unital. A \textit{regular multiplier Hopf ($\,^*$-)algebra} consists
of a couple $(A,\Delta)$, with $A$ a non-degenerate ($^*$-)algebra,
and $\Delta$, the \textit{comultiplication}, a morphism from $A$ to
$A\odot A$, where $\odot$ denotes the algebraic tensor product.
Moreover, $(A,\Delta)$ has to satisfy the following conditions:
\begin{itemize}\item[M.1] $(\Delta\otimes \iota)\Delta = (\iota\otimes
\Delta)\Delta$ \qquad (coassociativity). \item[M.2] The maps
\begin{itemize}
\item[]{ $T_{\Delta 2}: A \odot A \rightarrow M(A \odot A):
a\odot b \rightarrow \Delta (a)(1 \odot b)$, }

\item[]{ $T_{1 \Delta}: A \odot A \rightarrow M(A \odot A):
a\odot b \rightarrow (a \odot 1) \Delta (b)$, }

\item[]{ $T_{\Delta 1}: A \odot A \rightarrow M(A \odot A):
a\odot b \rightarrow \Delta (a) (b\odot 1)$, }
\item[]{$T_{2 \Delta}: A \odot A \rightarrow M(A \odot A):
a\odot b \rightarrow (1\odot a)\Delta (b)$ }
\end{itemize} all induce linear bijections $A\odot A\rightarrow A\odot A$.
\end{itemize}
Here, and elsewhere in the text, we will use $\iota$ to denote the
identity map.\\

\noindent The $T$-maps can be used to define a co-unit $\varepsilon$
(which will be a \linebreak ($^*$-)homomorphism from $A$ to
$\mathbb{C}$) and an antipode $S$ (which will be a linear
anti-morphism, satisfying $S(S(a)^*)^*=a$ for all $a\in A$ when $A$
carries an involution). Both co-unit and antipode will be unique,
and will satisfy the corresponding equations of those defining them
in the Hopf
algebra case.\\

\textit{($\,^*$-)Algebraic quantum groups}\\

\noindent Our second object forms an intermediate step between the
former, purely algebraic notion of a regular multiplier Hopf
algebra, and the analytic set-up of a locally compact quantum group.
A \textit{($\,^*$-)algebraic quantum group} is a regular multiplier
Hopf ($^*$-)algebra $(A,\Delta)$, for which there exists a non-zero
linear functional $\varphi$ on $A$, such that
\[(\iota\otimes \varphi)(\Delta(a)(b\otimes 1))=\varphi(a)b,\qquad
\textrm{for all } a,b\in A.\] When $A$ is a $^*$-algebra, we
\textit{impose the extra condition of positivity on $\varphi$:} for
every $a\in A$, we have $\varphi(a^*a)\geq 0$. This extra
condition is in fact very restrictive, as we shall see.\\

\noindent We can prove that $\varphi$ is unique up to multiplication
with a scalar. It will be faithful in the following sense: if
$\varphi(ab)=0$ for all $b\in A$, or $\varphi(ba)=0$ for all $b\in
A$, then $a=0$. The functional $\varphi$ is called the \emph{left
Haar functional}. Then $(A,\Delta)$ will also have a non-zero
functional $\psi$, again unique up to a scalar, such that
\[(\psi\otimes \iota)(\Delta(a)(1\otimes b))=\psi(a)b,
\qquad\textrm{for all } a,b\in A.\] We will call $\psi$ a
\emph{right Haar functional}. If $A$ is a $^*$-algebraic quantum
group, we can still choose $\psi$ to be positive. We note however,
that to arrive at this functional, a detour into analytic landscape
(with aid of the GNS-device for $\varphi$) seemed inevitable. The
problem is of course that the evident right Haar functional
$\psi=\varphi\circ S$ is not necessarily positive. To create the
right $\psi$, we needed some analytic machinery, namely the square
root of the modular element, or a polar decomposition of the
antipode (see \cite{Kus1}). In this paper we show, that it is
possible to arrive at the positivity of the right invariant
functional by \textit{purely algebraic means} (see the remark after
Theorem \ref{the2}). This means that $^*$-algebraic quantum groups
are appropriate objects of study for
algebraists with a fear of analysis.\\

\noindent Algebraic quantum groups have some really nice features.
For example, there exists a unique automorphism $\sigma$ of the
algebra $A$, satisfying $\varphi(ab)=\varphi(b\sigma(a))$ for all
$a,b\in A$. It is called \emph{the modular automorphism}, a notion
coming from the theory of weights on von Neumann-algebras (see
below). There also exists a unique multiplier $\delta$ such that
\[(\varphi \otimes \iota)(\Delta(a)(1\otimes b))=\varphi(a)\delta
b,\]\[(\varphi\otimes \iota)((1\otimes
b)\Delta(a))=\varphi(a)b\delta,\] for all $a,b\in A$. It is called
\emph{the modular element}, as it is the non-commutative equivalent
of the modular function in the theory of locally compact groups.
When $A$ is a $^*$-algebraic quantum group, $\delta$ will indeed be
a positive element (i.e. $\delta=q^*q$ for
some $q\in M(A)$).\\

\noindent There also is a particular number that can be associated
with an algebraic quantum group. Since $\varphi\circ S^2$ is a left
Haar functional, the unicity of $\varphi$ implies there exists
$\mu\in \mathbb{C}$ such that $\varphi(S^2(a))=\mu \varphi(a)$, for
all $a\in \mathbb{C}$. This number $\mu\in \mathbb{C}$ is called
\textit{the scaling constant} of $(A,\Delta)$. In an early stage,
examples of algebraic quantum groups were found where $\mu\neq 1$
(see \cite{VDae2}). However, it remained an open question whether
$^*$-algebraic quantum groups existed with $\mu\neq 1$. We will show
in this paper that in fact $\mu=1$ for all $^*$-algebraic quantum
groups (see Theorem \ref{the1}).
\\

\textit{C$^*$-algebraic quantum groups}\\

\noindent A \textit{reduced C$^*$-algebraic quantum group} is the
non-commutative version of the space of continuous complex
functions, vanishing at infinity, of a locally compact group. It
consists of a couple $(A,\Delta)$, with $A$ a C$^*$-algebra,
encoding the topology of the quantum space, and $\Delta$ a morphism
$A\rightarrow A\otimes A$, encoding the group-multiplication in the
quantum space. Here $\otimes$ denotes the minimal C$^*$-algebraic
tensor product. Moreover, $(A,\Delta)$ has to satisfy the following
conditions, which are the appropriate non-commutative translations
of the associativity condition and the cancelation
property\footnote{Actually, the cancelation property says that
$\Delta(A)(1\otimes A)$ and $\Delta(A)(A\otimes 1)$ are dense in
$A\otimes A$, but it can be shown that the given (weaker) condition
implies this.} of a locally compact group:
\begin{itemize}\item[C.1] $(\Delta\otimes \iota)\Delta = (\iota\otimes
\Delta)\Delta$\qquad (coassociativity).\item[C.2] $\{(\omega\otimes
\iota)\Delta(a)\mid a\in A,\,\omega\in A^*\}$ and $\{(\iota\otimes
\omega)\Delta(a)\mid a\in A,\,\omega\in A^*\}$ are norm-dense in
$A$.\end{itemize} Here $A^*$ denotes the space of continuous
functionals on $A$.\\

\noindent In contrast with the theory of locally compact groups, it
seems unlikely that a satisfactory theory of locally compact quantum
groups can be built from simply topological and algebraic axioms.
Some analytic structure has to be included. More precisely, we have
to assume the existence of invariant weights, which correspond to
the notion of Haar measures on a locally compact group. The final
axiom reads:
\begin{itemize}\item[C.3] There exist faithful KMS-weights $\varphi$ and
$\psi$ on $A$, such that
\begin{enumerate} \item $\varphi((\omega\otimes
\iota)\Delta(a))=\varphi(a)\omega(1)$ for $a\in A^+$ and $\omega\in
A_+^*$, \item $\psi((\iota\otimes
\omega)\Delta(a))=\psi(a)\omega(1)$ for $a\in A^+$ and $\omega \in
A_+^*$.\end{enumerate}\end{itemize} Actually, this axiom can be
relaxed in a non-trivial manner (see \cite{Kus2}). For an account of
the theory of weights
on C$^*$-algebras, we refer the reader to \cite{Kus4}.\\

\noindent It turns out that $\varphi$ and $\psi$ are unique up to
multiplication with a positive scalar. They are called respectively
a \textit{left} and \textit{right Haar weight.} Unlike the case of
algebraic quantum groups, both left and right Haar weight must be
part of the definition to have a good theory. The most important
objects associated with $(A,\Delta)$ are the \emph{multiplicative
unitaries}, which essentially carry all information about
$(A,\Delta)$. To introduce them, we first recall the notion of the
GNS-representation associated to $\psi$. Denote the algebra of
$\psi$-integrable elements by $\mathscr{M}_\psi$ and the set of
square integrable elements by $\mathscr{N}_\psi$. The GNS-space
associated to $\psi$ is the closure $\mathscr{H}_\psi$ of the
pre-Hilbert space $\mathscr{N}_\psi$, with scalar product defined by
$\langle a,b\rangle = \psi(b^*a)$ for $a,b\in \mathscr{N}_\psi$. The
injection of $\mathscr{N}_\psi$ into $\mathscr{H}_\psi$ will be
denoted by $\Lambda_\psi$. We can construct a faithful
representation of $A$ on $\mathscr{H}_\psi$ via left multiplication.
We can do the same for $\varphi$, obtaining a Hilbert space
$\mathscr{H}_\varphi$ and an injection $\Lambda_\varphi:
\mathscr{N}_\varphi\rightarrow \mathscr{H}_\varphi$. Moreover, there
exists a unitary from $\mathscr{H}_\psi$ to $\mathscr{H}_\varphi$,
commuting with the action of $A$. Since the representations are
faithful, we can identify both Hilbert spaces, and denote this
Hilbert space as $\mathscr{H}$. We will let elements of $A$ act
directly on $\mathscr{H}$ as operators (suppressing the
representation). Now we can define the multiplicative unitary $W$,
also called the \emph{left regular representation}: it is the
unitary operator on $\mathscr{H}\otimes\mathscr{H}$, characterized
by
\[(\iota\otimes
\omega)(W^*)\Lambda_\varphi(x)=\Lambda_\varphi((\iota\otimes
\omega)\Delta(x)),\qquad x\in \mathscr{N}_\varphi \textrm{ and }
\omega \in B(\mathscr{H})_*.\] Here $B(\mathscr{H})$ denotes the
space of bounded operators on $\mathscr{H}$ and $B(\mathscr{H})_*$
its pre-dual. It implements the comultiplication as follows:
\[\Delta(x)=W^*(1\otimes x)W, \qquad\textrm{for all } x\in M.\]
Moreover, the normclosure of the set
$\{(\iota\otimes \omega)(W)\mid \omega\in B(\mathscr{H})_*\}$ will be equal to $A$.\\

\noindent We can also define the multiplicative unitary $V$, called
\textit{the right regular representation}: it is the unitary
operator on $\mathscr{H}\otimes\mathscr{H}$, determined by
\[(\omega\otimes \iota)(V)\Lambda_\psi(x)=\Lambda_\psi((\omega\otimes
\iota)\Delta(x)),\qquad x\in \mathscr{N}_\psi \textrm{ and } \omega
\in B(\mathscr{H})_*.\] It also implements the comultiplication:
\[\Delta(x)=V(x\otimes 1)V^*, \qquad\textrm{for all } x\in M.\]
Again, the normclosure of the set $\{(\omega\otimes \iota)(V)\mid
\omega\in B(\mathscr{H})_*\}$ will be equal to $A$.\\

\noindent With the aid of a multiplicative unitary, it is possible
to construct a dual quantum group $(\hat{A},\hat{\Delta})$. Just as
in the classical case, the bi-dual will be isomorphic to the
original quantum group.\\

\noindent The right regular representation can also be used to
define the \textit{antipode} on $(A,\Delta)$. It is the (possibly
unbounded) closed linear map $S$ from $A$ to $A$, with a core
consisting of elements of the form $(\omega\otimes \iota)(V)$,
$\omega\in B(\mathscr{H})_*$, such that \[S((\omega\otimes
\iota)(V))=(\omega\otimes \iota)(V^*).\] This map has a polar
decomposition, consisting of a (point-wise) \textit{normcontinuous}
one-parameter group $\tau$ on $A$ (called the \emph{scaling group})
and a $^*$-anti-automorphism $R$ of $A$ (called the \emph{unitary
antipode}). Then the antipode equals the map $R\circ \tau_{-i/2}$,
where $\tau_{-i/2}$ is the analytic continuation of $\tau$ to the
point $-i/2$.\\

\noindent Since $\psi$ and $\varphi$ are KMS-weights, there exist
respective normcontinuous one-parameter groups $\sigma'$ and
$\sigma$ on $A$, called the \emph{modular one-parameter groups}
(associated with $\psi$ and $\varphi$). It can be shown that there
exists a (possibly unbounded) positive operator $\delta$ on
$\mathscr{H}$, affiliated with $A$ (in the sense of Woronowicz),
such that $\psi=\varphi(\delta^{1/2}\cdot\delta^{1/2})$ for all
$a\in A$ and $t\in \mathbb{R}$. For a correct interpretation of this
equality, we refer to \cite{Vae1}. We call $\delta$ the
\textit{modular element} of $(A,\Delta)$. In the sequel, we will
frequently use the commutation rules between these objects and
$\Delta$:

\[\left.\begin{array}{llll} \Delta \tau_t = (\tau_t\otimes
\tau_t)\Delta &&& \Delta \sigma_t = (\tau_t\otimes \sigma_t)\Delta
\\ \Delta\sigma'_t = (\sigma'_t\otimes \tau_{-t})\Delta &&& \Delta \tau_t =
(\sigma_t\otimes \sigma'_{-t})\Delta. \end{array} \right.\]\\[4pt]

\noindent There also exists $\nu \in \mathbb{R}^+$ such that
$\sigma_t(\delta)=\nu^{t}\delta$. This constant $\nu$ is called the
\emph{scaling constant}. It arises naturally in the framework
considered by Kustermans and Vaes. Therefore, it was an important
question whether there exist locally compact quantum groups where
this constant is not trivially 1. Such quantum groups do indeed
exist: an interesting example is the quantum $az+b$-group (see
\cite{VDae7}).\\

\textit{Von Neumann-algebraic quantum groups}\\

\noindent Our final object is a \textit{von Neumann-algebraic
quantum group}. It consists of a couple $(M,\Delta)$, with $M$ a von
Neumann-algebra, and $\Delta$ a normal morphism from $M$ to
$M\otimes M$, satisfying the coassociativity condition. Here
$\otimes$ denotes the ordinary von Neumann-algebraic tensor product.
We also assume that there exist faithful, semi-finite normal weights
$\varphi$ and $\psi$ on $M$, satisfying
\begin{enumerate} \item $\varphi((\omega\otimes
\iota)\Delta(a))=\varphi(a)\omega(1)$ for $a\in M^+$ and $\omega\in
M_{*,+}$, \item $\psi((\iota\otimes
\omega)\Delta(a))=\psi(a)\omega(1)$ for $a\in M^+$ and $\omega \in
M_{*,+}$.\end{enumerate} The theory then develops parallel to the
C$^*$-algebraic case, but with some significant simplifications (see
\cite{VDae4}). Remark for example that its definition does not
mention the `cancelation law': this law will automatically be
fulfilled! However, it can be shown that it contains as much
information as a C$^*$-algebraic locally compact quantum group:
there is a canonical bijection between C$^*$-algebraic and von
Neumann-algebraic quantum groups. When working in the von
Neumann-algebraic context, we will use the same notation as in the
C$^*$-algebraic case (the antipode will be denoted by $S$, the
scaling group by $\tau$,
...). For the theory of weights on von Neumann-algebras, we refer the reader to \cite{Tak1}.\\

\section{Multiplier Hopf algebras imbedded in locally compact quantum groups}

In this section, we fix a C$^*$-algebraic quantum group $(A,\Delta)$
and a regular multiplier Hopf algebra $(A_0,\Delta_0)$. The von
Neumann-algebraic quantum group associated to $(A,\Delta)$ will be
denoted by $(M,\Delta)$. We will use notations as before, but the
structural maps for $A_0$ will be indexed by 0 (whenever this causes
no confusion). We also fix a left Haar weight $\varphi$ on
$(A,\Delta)$. As a right Haar weight on $(A,\Delta)$, we choose
$\psi=\varphi\circ R$. When using the notation
$\mathscr{N}_\varphi$, we will always specify whether we mean the
square integrable elements in $M$ or in $A$.
\\[4pt]

\qquad \Ass $A_0\subseteq M$.\\

\noindent This means that $A_0$ is a subalgebra of $M$, not
necessarily invariant under the $^*$-involution. We also want to
impose a certain compatibility between $\Delta$ and $\Delta_0$, but
we have to be careful: $M(A_0)$ bears no natural relation to $M$.
For example, denoting by $j$ the inclusion of $A_0$ in $M$, the
identity $(j\otimes j)\circ \Delta_0=\Delta\circ j$ can be
meaningless if $j$ has no well-defined extension to $M(A_0)$. We
will however assume the following: for all $a,b\in A_0$,
\[\Delta_0(a)(1\otimes
b)=\Delta(a)(1\otimes b),\]\[\Delta_0(a)(b\otimes
1)=\Delta(a)(b\otimes
1),\]\[(a\otimes1)\Delta_0(b)=(a\otimes1)\Delta(b),\]\[(1\otimes
a)\Delta_0(b)=(1\otimes a)\Delta(b).\]

\noindent \textit{Remark.} This condition is strictly weaker than
the condition $(j\otimes j)\circ \Delta_0=\Delta\circ j$, when it
makes sense. For example, the imbedding of $C(\mathbb{Z})$ in
$C(\mathbb{Z})$ sending $\delta_i$ to $\delta_{2i}$ satisfies the
former, but not the latter
condition.\\

\noindent Our first result shows that the antipode $S$ of $M$
restricts to the antipode $S_0$ of $A_0$. The hard part consists of
showing that $A_0$ lies in the domain of $S$. We will need a lemma
which is interesting in its own right. It is a kind of cancelation
property involving $M$ and $\hat{M}'$, the commutant of the dual
quantum group $\hat{M}$.

\noindent
\begin{Lemma}\label{lem3} Suppose $a\in M$ and $x\in \hat{M}'$ satisfy $ax=0$. Then $a=0$ or
$x=0$.

\begin{proof}

Let $W$ be the left regular representation for $M$. We recall that
$W\in M\otimes \hat{M}$ (see e.g. \cite{Kus3}). So if $xa=0$, then
\begin{eqnarray*} W^*(1\otimes xa)W &=& (1\otimes x)W^*(1\otimes
a)W\\&=& (1\otimes x)\Delta(a)\\ &=&0.\end{eqnarray*} Assume $x\neq
0$. Choose $\omega\in B(\mathscr{H})_{*,+}\,$ such that
$\omega(xx^*)=1$. Then we have $(\iota \otimes \omega(x\cdot
x^*))\Delta(aa^*)=0$. Applying $\psi$ and using the strong right
invariance property, we get $\psi(aa^*)\omega(xx^*)=\psi(aa^*)=0$.
Since $\psi$ is faithful, $a$ must be zero.

\end{proof}
\end{Lemma}

\noindent \textit{Remark.} Applying $\hat{J}\cdot\hat{J}$, we see
that also the following is true: if $a\in M$ and $x\in \hat{M}$,
then $ax=0$
implies either $a=0$ or $x=0$.\\

\noindent We can show now that the antipodes of $M$ and $A_0$
coincide.

\begin{Proposition} $A_0$ lies in the domain of $S$, and $S_{|A_0}$ will
be the antipode of $(A_0,\Delta_0)$.

\begin{proof}

Let $b$ be an element of $A_0$. We will show that $b\in
\mathscr{D}(S)$ and $S(b)=S_0(b)$. We start by choosing some fixed
$a$ in $A_0$. We can pick $p_i,q_i$ in $A_0$ such that
\[a\otimes b = \sum_{i=1}^n (p_i\otimes 1)\Delta(q_i).\] Then
\[\Delta(a)(1\otimes S_0(b)) = \sum_{i=1}^n \Delta(p_i)(q_i\otimes
1).\]\\

\noindent Let $y$ be $(\omega_{c,d}(a\cdot) \otimes \iota)(V) =
(\psi\otimes \iota)((c^*a\otimes 1)\Delta(d))$, where $c,d$ are
square integrable elements in $M$ and $\omega_{c,d}=\langle\cdot
\Lambda_{\psi}(d),\Lambda_{\psi}(c)\rangle$. Then
\begin{eqnarray*} by&=&(\psi\otimes \iota)((c^*a\otimes b) \Delta(d))\\&=&(\psi\otimes\iota)\sum (c^*p_i\otimes
1)\Delta(q_id).\end{eqnarray*} We know that this last expression is
in $\mathscr{M}_{\psi\otimes \iota}$ (where we use the
slice-notation) and that
\[(\psi\otimes \iota)\sum (c^*p_i\otimes 1)\Delta(q_id)\in
\mathscr{D}(S),\] with \[ S((\psi\otimes\iota)\sum(c^*p_i\otimes
1)\Delta(q_id)) =(\psi\otimes \iota)\sum
\Delta(c^*p_i)(q_id\otimes1).\] So
\begin{eqnarray*} S(by)&=&
(\psi\otimes \iota)\sum \Delta(c^*p_i)(q_id\otimes1)
\\&=&(\psi\otimes \iota)( \Delta(c^*a)(d\otimes
S_0(b)))\\&=&S(y)S_0(b).\end{eqnarray*}\\

\noindent Denote by $C$ the linear span of all such $y$, with $c$
and $d$ varying. We show that $C$ is an ultra-strong$^*$-core for
$S$. First remark that functionals of the form
$\omega_{c,d}(a\cdot)$ have a norm-dense linear span in
$(\hat{M}')_*$. Indeed: if $z\in\hat{M}'$ such that $\langle az
\Lambda_{\psi}(d),\Lambda_{\psi}(c)\rangle=0$ for all $c,d\in
\mathscr{N}_\psi$, then $az=0$, hence $z=0$ by the previous lemma.
Then, since $V\in \hat{M}'\otimes M$, there exists for every
$\omega\in B(\mathscr{H})_*$ a sequence of $c_n,d_n\in
\mathscr{N}_\psi$ such that $(\omega_{c_n,d_n}(a\cdot)\otimes
\iota)(V)\rightarrow (\omega\otimes \iota)(V)$ and
$(\omega_{c_n,d_n}(a\cdot)\otimes \iota)(V^*)\rightarrow
(\omega\otimes \iota)(V^*)$. Since $\{(\omega\otimes \iota)(V)\mid
\omega\in B(\mathscr{H})_*\}$ is an ultra-strong$^*$-core for $S$,
the same
will be true for $C$.\\

By choosing a net $y_\alpha$ in $C$ such that $y_\alpha\rightarrow
1$ and $S(y_\alpha)\rightarrow 1$ in the ultra-strong$^*$-topology,
we can conclude that $b\in \mathscr{D}(S)$ and $S(b)=S_0(b)$.

\end{proof}\end{Proposition}

\noindent The previous proposition implies that $A_0\subset
\mathscr{D}(\tau_{z})$ for every $z$ in $\mathbb{C}$, i.e. every
$a\in A_0$ is analytic with respect to $\tau$. Indeed: $a\in
\mathscr{D}(S)$ means that $a\in \mathscr{D}(\tau_{-i/2})$. Since
$S(S_0^{-1}(a))=a$ for $a\in A_0$, we also have that $a\in
\mathscr{D}(S^{-1})=\mathscr{D}(\tau_{i/2})$. So $A_0\subset
\mathscr{D}(\tau_{ni})$ for every integer $n\in \mathbb{Z}$.\\

\noindent This again illustrates the lack of analytic structure of a
general algebraic quantum group: if its antipode $S$ satisfies
$S^{2n}=\iota$, but $S^2\neq\iota$, then it can not be imbedded in a
locally compact quantum group \emph{at all}. It can not be \emph{at
all}. Such algebraic quantum groups do indeed exist (see e.g.
\cite{VDae2}).\\[4pt]

\noindent We can also use Lemma \ref{lem3} to prove that actually
$A_0\subset M(A)$. Fix $a\in A_0$. Choose $b\in A_0$ and $\omega\in
B(\mathscr{H})_*$.
Then \begin{eqnarray*} a\otimes b&=&\sum (q_i\otimes 1)\Delta(p_i)\\
&=& \sum (q_i\otimes 1)V(p_i\otimes 1)V^*,\end{eqnarray*} for some
$p_i,q_i$ in $A_0$. Multiplying from the right with $V$ and applying
$\omega\otimes \iota$, we get $b(\omega(a\cdot)\otimes \iota)(V)\in
A$. But as we have shown, the set $\{(\omega(a\cdot)\otimes
\iota)(V)\mid \omega\in B(\mathscr{H})_*\}$ is norm-dense in $A$.
Hence $bA\subseteq A$.
Similarly $Ab\subseteq A$, and thus $A_0\subset M(A)$.\\

\noindent As a second important result, we show that $A_0$ consists
of analytic elements for $\sigma$. This follows easily from the
following proposition, which elucidates the behavior of $A_0$ with
respect to the one-parameter group $\kappa$, where
$\kappa_t=\sigma_t\tau_{-t}$. It will be decisive in obtaining some
structural properties of $^*$-algebraic quantum groups, as we will
show in the third section.

\begin{Proposition}\label{lem1} $A_0\subset\mathscr{D}(\kappa_z)$
for all $z\in\mathbb{C}$, and $\kappa_z(A_0)\subset A_0$. Here
$\kappa_z$ denotes the analytic continuation of the one-parameter
group $\kappa$ to the point $z\in \mathbb{C}$, and
$\mathscr{D}(\kappa_z)$ denotes its domain.

\begin{proof}

Let $b$ be a fixed element of $A_0$.  Choose a non-zero $a\in A_0$,
and write \[a\otimes b = \sum_{i=1}^n \Delta(p_i)(1\otimes q_i), \]
with $p_i,q_i\in A_0$. Using the commutation relations between
$\Delta$, $\tau$, $\sigma$ and $\sigma'$, we get that
\[\kappa_{-t}(a)\otimes \rho_t(b)=\sum \Delta(p_i)(1\otimes \rho_t(q_i)), \qquad\textrm{for all } t\in
\mathbb{R},\] where $\rho_t=\sigma'_t\tau_t$. Choose $c\in A_0$ such
that $cb\neq 0$, and multiply this equation to the left with
$1\otimes c$ to get
\[\kappa_{-t}(a)\otimes c\rho_t(b)=\sum
((1\otimes c)\Delta(p_i))(1\otimes \rho_t(q_i)).\] Choose
$a_{ij},b_{ij}\in A_0$ such that \[(1\otimes
c)\Delta(p_i)=\sum_{j=1}^{m_i} a_{ij}\otimes b_{ij},\] and let $L$
be the finite-dimensional space spanned by the $a_{ij}$. We see that
$\kappa_{-t}(a)\otimes c\rho_t(b)\in L\odot M$, for every $t\in
\mathbb{R}$. Since $c\rho_0(b)=cb\neq 0$, and $\rho_t$ is strongly
continuous, we get that there exists a $\delta>0$ such that
$c\rho_t(b)\neq 0$ for all $t$ with $|t|<\delta$. This means
$\kappa_t(a)\in L$ for all $|t|<\delta$.

For every $\varepsilon>0$, let
$K_\varepsilon=\textrm{span}\{\kappa_t(a)\mid |t|<\varepsilon\}$,
and $n_\varepsilon=\textrm{dim}(K_\varepsilon)$. For small
$\varepsilon$, we have $n_\varepsilon\in \mathbb{N}$. Choose an
$\varepsilon$ where this dimension reaches a minimum. Then
$K=K_{\varepsilon}=K_{\varepsilon/2}$ will be a finite-dimensional
subspace containing $a$, invariant under $\kappa_t$, for all
$t\in\mathbb{R}$.

Now $\kappa$ induces a continuous homomorphism
$\tilde{\kappa}:\mathbb{R}\rightarrow \textrm{GL}(K)$. It is a
well-known fact that such a homomorphism is necessarily analytic.
Thus $a\in \mathscr{D}(\kappa_z)$, and $\kappa_z(a)\in K\subseteq
A_0$. This concludes the proof.

\end{proof}
\end{Proposition}

\noindent \textit{Remarks.} (i)\phantom a We can actually show that
$A_0$ is spanned by eigenvectors for $\kappa_t$. For every continuous one-parameter group of isometries on a finite-dimensional Banach space is diagonizable.
Namely, using
notation as in the proof, consider a Banach space isomorphism of $K$
to $\mathbb{C}^n$ with the usual Hilbert space structure. Then
$\kappa_t$ will be transformed to a uniformly bounded,
norm-continuous one-parameter group of operators on $\mathbb{C}^n$.
Since $\mathbb{R}$ is amenable, we can choose another Hilbert
structure on $\mathbb{C}^n$ such that the one-parameter group will
consist of unitaries. Hence the action is diagonizable.\\[2pt]
\indent\indent\phantom a\phantom a (ii)\phantom a The lemma remains
true if we replace $\kappa_t$ by
$\rho_t=\tau_t\sigma'_t$ or $\sigma_t\sigma'_t$.\\

\begin{Cor}\label{lem2}
$A_0$ consists of analytic elements for $\sigma$.
\begin{proof}

This follows easily from the previous two statements. If $a\in A_0$,
we know that $a$ is analytic for $\tau_t$ and
$\kappa_t=\sigma_t\tau_{-t}$. If $z\in \mathbb{C}$, then
$\tau_z(\kappa_z(a))$ makes sense, since $A_0$ is invariant under
$\kappa_z$. Since $\sigma_z$ is the closure of $\tau_z\circ
\kappa_z$, we arrive at $a\in \mathscr{D}(\sigma_z)$.

\end{proof}
\end{Cor}

\noindent As a consequence, $A_0$ is invariant under $\sigma_{ni}$
and
$\sigma'_{ni}$, with $n\in \mathbb{Z}$.\\

\noindent \textit{Remark.} We do not know if $A_0$, or even the von
Neumann-algebra $N$ generated by it, has to be invariant under the
one-parameter groups $\sigma$ and $\tau$. There seems to be an
analytic obstruction to be able to conclude this. It is however easy
to see that if $N$ is invariant under either $\sigma$, $\tau$ or
$\delta^{-it}\cdot \delta^{it}$, then it is invariant under all of
them (see e.g. Proposition \ref{prop5}).\\

\noindent Next, we impose a stronger condition on $A_0$:\\

\qquad \Ass $A_0\subseteq A$.\\

\noindent We will say then that $A_0$ has a proper imbedding in $A$.
Because $A_0$ is now a subspace of the C$^*$-algebra $A$, we can say
more about its connection to $\varphi$. We first need a simple
lemma, which also appears in some form in \cite{Lan1}:

\begin{Lemma} Suppose that $a\in A\cap \mathscr{D}(\sigma_{i/2})$ and $e\in A$
satisfy $ea=a$. Then $a\in \mathscr{N}_\varphi$.
\end{Lemma}

\begin{proof}
Choose $c$ in $A\cap\mathscr{M}_\varphi^+$ such that
$0\leq\|c-e^*e\|\leq1/2$. This is possible because
$\mathscr{M}_\varphi^+\cap A$ is \textit{normdense} in $A^+$. Then
\begin{eqnarray*} \frac{1}{2} a^*a &\leq& a^*(1+c-e^*e)a \\ &=&
a^*ca.\end{eqnarray*}%
Since $a$ is analytic with respect to $\sigma$, we conclude that
$a^*ca\in \mathscr{M}_\varphi^+$, and thus $a^*a\in
\mathscr{M}_\varphi^+$.
\end{proof}

\begin{Proposition} $A_0$ belongs to the Tomita algebra of
$\varphi$: $A_0\subset \mathscr{T}_\varphi=\{x\in
\mathscr{N}_\varphi\cap \mathscr{N}_\varphi^*\mid x\in
\mathscr{D}(\sigma_z)\textrm{ and } \sigma(z)\in
\mathscr{N}_\varphi\cap\mathscr{N}_\varphi^*\textrm{ for all }z\in
\mathbb{C}\}$.
\begin{proof}
We know that $A_0$ has local units: for every $a\in A_0$ there exist
$e,f\in A_0$ such that $a=ea$ and $a=af$ (see e.g. \cite{VDae8}).
So, since $A_0$ consists of analytic elements for $(\sigma_t)$, we
can apply the previous lemma to each element of $(\bigcup_{z\in
\mathbb{C}} \sigma_z(A_0))$ and $A_0^*$. This implies that
$A_0\subset \mathscr{T}_\varphi$.
\end{proof}
\end{Proposition}

\noindent \textit{Remark.} A converse is also true. Suppose $A_0$
consists of square integrable elements in $M$. Then $A_0$ will be a
subset of $A$. Namely: let $b$ be a fixed element in $A_0$ such that
$\varphi(b^*b)=1$. Choose $a$ in $A_0$. Then $a\otimes b=\sum
\Delta(q_i)(p_i\otimes 1)$ with $p_i,q_i\in A_0$. Multiply to the
left with $1\otimes b^*$ and apply $\iota\otimes\varphi$, then
$a=\sum (\iota\otimes \langle\cdot
\Lambda_\varphi(q_i),\Lambda_\varphi(b)\rangle)(W^*)p_i$. Since
$A_0\subset M(A)$, this is an element of $A$.\\

\noindent The previous proposition has the interesting corollary
that the scaling constant of $A$ is necessarily trivial. We will
come back to this fact in the third section, where we apply our
techniques to $^*$-algebraic quantum groups (see Theorem
\ref{the1}).

\begin{Cor} The scaling constant $\nu$ of $(A,\Delta)$ equals
1.
\end{Cor}

\begin{proof}
We have that $\nu^{-\frac{1}{2}t}\kappa_t$ induces a one-parameter
unitary group $u_t$ on $\mathscr{H}$, where
$\kappa_t=\sigma_t\circ\tau_{-t}$. As in the proof of lemma
\ref{lem1}, there is a finite-dimensional subspace $K$ of $A_0$ that
is invariant under $\kappa$. Therefore $V=\Lambda_\varphi(K)$ is
invariant under $u$. This means that there exists a non-zero
$v=\Lambda_\varphi(x)\in V$ such that $u_t(v)=e^{it\lambda}v$, for
some $\lambda\in \mathbb{R}$. Hence $\nu^{-\frac{1}{2}t}\kappa_t(x)=
e^{it\lambda}x$. But, since $\kappa_t$ is a one-parameter group of
$^*$-automorphisms, we get
\begin{eqnarray*}
\|x\|&=&\|\kappa_t(x)\|\\&=&\|e^{it\lambda}\nu^{\frac{1}{2}t}x\|\\&=&\nu^{\frac{1}{2}t}\|x\|.\end{eqnarray*}
So $\nu=1$.\end{proof}¨\\

\noindent From the previous proposition, it follows that
$A_0\subset\mathscr{N}_\varphi\cap \mathscr{N}_\varphi^*$. Because
$A_0^2=A_0$, we also have $A_0\subset \mathscr{M}_\varphi$, so every
element of $A_0$ is integrable with respect to $\varphi$. However,
we can not conclude that $(A_0,\Delta_0)$ is an algebraic quantum
group, because we do not know if the restriction $\varphi_0$ of
$\varphi$ to $A_0$ is non-zero. In any case, it will be left
invariant: If $a,b\in A_0$, then $\Delta(a)(b\otimes 1)\in
\mathscr{M}_{\iota\otimes \varphi}$, and
\begin{eqnarray*} (\iota\otimes\varphi_0)(\Delta_0(a)(b\otimes
1))&=&(\iota\otimes \varphi)(\Delta(a)(b\otimes 1))\\ &=&
\varphi(a)b\\ &=& \varphi_0(a)b.\end{eqnarray*}\\[4pt]

\quad \begin{tabular}{ll}\Ass & $A_0\subseteq A$ and
$\varphi_{|A_0}\neq 0$.

\end{tabular}\\

\noindent The assumption is sufficient to conclude that
$(A_0,\Delta_0)$ is an algebraic quantum group, as we have shown.
Remark that the second condition is automatically fulfilled if
$(A_0,\Delta_0)$ is a multiplier Hopf $^*$-algebra (with the same
$^*$-involution as in $A$).\\

\noindent We now show that $A_0$ itself possesses an analytic
structure, thus generalizing the results in \cite{Kus5}.

\begin{Proposition}  Let $\delta$ be the modular element of $(A,\Delta)$, and $\delta_0$
the modular element for $A_0$. Then every $a$ in $A_0$ is a left and
a right multiplier for $\delta$ such that $a\delta=a\delta_0$ and
$\delta a=\delta_0 a$. Moreover, we have that $\delta^{z}A_0= A_0$
and $A_0\delta^{z}= A_0$.

\begin{proof} Choose a fixed $b$ in $A_0$ with
$\varphi(b)\neq0$. Choose $a$ in $A_0$. Then $a\otimes b=\sum
\Delta(p_i)(q_i\otimes 1)$ for certain $p_i,q_i$ in $A_0$.
Multiplying to the left with $\delta^{it}\otimes \delta^{it}$, we
get that
\[\delta^{it}a\otimes\delta^{it}b=\sum\Delta(\delta^{it}p_i)(q_i\otimes
1).\] Every term in the right hand side lies in
$\mathscr{M}_{\iota\otimes \varphi}$, so applying $\iota\otimes
\varphi$ to each side, we get
\[\varphi(\delta^{it}b)\delta^{it}a=\sum\varphi(\delta^{it}p_i)q_i\in
A_0.\]\\

Denote by $L$ the finite-dimensional vector space spanned by the
$q_i$. Since $t\rightarrow \varphi(\delta^{it}b)$ is a continuous
function and $\varphi(b)=1$, we can choose $t$ small such that
$\varphi(\delta^{it}b)\neq 0$. For such $t$ we have $\delta^{it}a\in
L$. A similar argument as in Proposition \ref{lem1} let's us
conclude that the linear span of the $\delta^{it}a$, with
$t\in\mathbb{R}$, is a finite-dimensional subspace of $A_0$. This
easily implies that $a$ is a right multiplier of $\delta$ and
$\delta a\in A_0$. So every element of $A_0$ lies in the domain of
left multiplication with $\delta^z$, and $\delta^{z}A_0=A_0$. Since
the space $B_0=\{a\in A\mid a^*\in A_0\}$ also has the structure of
a multiplier Hopf algebra, and $\varphi_{|B_0}\neq 0$, we can
conclude that $B_0$ consists of right multipliers for $\delta$. So
$A_0$ consists of left multipliers for $\delta$, and
$A_0\delta^{z}=A_0$.\\

Choose a fixed $a\in A_0$ with $\varphi(a)\neq 0$. Let $b,c$ be
elements in $A_0$. We know that (see \cite{Kus2})
\begin{eqnarray*} \varphi((\iota\otimes
\langle\cdot\Lambda_\varphi(b),\Lambda_\varphi(c^*)\rangle
)\Delta(a))&=& \varphi(a)\langle
\delta^{1/2}\Lambda_\varphi(b),\delta^{1/2}\Lambda_\varphi(c^*)\rangle\\
&=& \varphi(a)\langle \Lambda_\varphi(b),\Lambda_\varphi(\delta c^*)\rangle\\
&=& \varphi_0(a)\varphi_0(c\delta b).\end{eqnarray*} On the other
hand,
\begin{eqnarray*} \varphi((\iota\otimes
\langle\cdot\Lambda_\varphi(b),\Lambda_\varphi(c^*)\rangle
)\Delta(a))&=& (\varphi\otimes \varphi)((1\otimes
c)\Delta(a)(1\otimes b))\\ &=& (\varphi_0\otimes
\varphi_0)((1\otimes c)\Delta_0(a)(1\otimes b))\\ &=&
\varphi_0(a)\varphi_0(c\delta_0b).\end{eqnarray*} Since $\varphi_0$
is faithful, $\delta_0 b=\delta b$ and $b\delta_0=b\delta$ for all
$b\in A_0$.

\end{proof}
\end{Proposition}

\noindent \textit{Remark.} In general (i.e. when $A_0\subset M(A)$),
we do not have to expect nice behavior of $A_0$ with respect to
$\delta$. Consider
for example the trivial quantum group $\mathbb{C}1$ in $M(A)$.\\

\noindent As we have remarked, the invariance under the
one-parameter groups of $A_0$ follows easily.

\begin{Proposition}\label{prop5}
$\tau_z(A_0)=\sigma_z(A_0)=R(A_0)=A_0$, for all $z\in \mathbb{C}$.
\begin{proof}

We have
$\sigma_{2z}(a)=\delta^{-iz}(\sigma_z'\sigma_z(a))\delta^{iz}$. But
$\sigma_z'\sigma_z$ and $\delta^{-iz}\cdot \delta^{iz}$ leave $A_0$
invariant. Hence $\sigma_z(A_0)\subseteq A_0$. Then also
$\tau_z=(\tau_z\sigma_{-z})\circ \sigma_z$ leaves $A_0$ invariant.
Since $R=S\circ \tau_{i/2}$, we have that $R$ leaves $A_0$
invariant.

\end{proof}
\end{Proposition}

\noindent Gathering all we have proven so far, we obtain the
following theorem:

\begin{Theorem} Let $(A,\Delta)$ be a reduced C$^*$-algebraic quantum group
with left Haar weight $\varphi$. Let $(A_0,\Delta_0)$ be a regular
multiplier Hopf algebra in $A$, such that \[\Delta_0(a)(1\otimes
b)=\Delta(a)(1\otimes b),\]\[\Delta_0(a)(b\otimes
1)=\Delta(a)(b\otimes
1),\]\[(a\otimes1)\Delta_0(b)=(a\otimes1)\Delta(b),\]\[(1\otimes
a)\Delta_0(b)=(1\otimes a)\Delta(b),\] for all $a,b\in A_0$. Then
$A_0$ will consist of integrable elements for $\varphi$. If
$\varphi_{|A_0}\neq0$, then $(A_0,\Delta_0)$ will be an algebraic
quantum group with left Haar functional $\varphi_0=\varphi_{|A_0}$.
Moreover, $A_0$ will consist of analytic elements for the modular
automorphism group, the scaling group, the unitary antipode and left
and right multiplication with the modular element of $(A,\Delta)$,
and $A_0$ will be invariant under all these actions.\end{Theorem}

As a corollary, we have

\begin{Cor} Let $(A,\Delta)$ be a reduced C$^*$-algebraic quantum group with
a dense, properly imbedded regular multiplier Hopf $^*$-algebra
$(A_0,\Delta_0)$. Then $(A_0,\Delta_0)$ is a $^*$-algebraic quantum
group, with associated C$^*$-algebraic quantum group
$(A,\Delta)$.\end{Cor}

\begin{proof} From the foregoing, we know that $(A_0,\Delta_0)$ is a
$^*$-algebraic quantum group with left Haar functional
$\varphi_0=\varphi_{|A_0}$. The only difficult step left to show, is
that $A_0$ is actually a core for the GNS-map $\Lambda_\varphi$. The
proof of this follows along the lines of Theorem 6.12. of
\cite{Kus1}.\\

Let $\Lambda_0$ be the closure of the restriction of
$\Lambda_\varphi$ to $A_0$. Choose a bounded net $(e_j)$ in $A_0$
converging strictly to 1. We can replace $e_j$ by
$\frac{1}{\sqrt{\pi}}\int \textrm{exp}(-t^2)\sigma_t(e_j)dt$, since
each will be an element of $A_0$ (because $\{ \sigma_t(e_j)\mid t\in
\mathbb{R}\}$ only spans a finite-dimensional space in $A_0$), and
the net will still be bounded, converging strictly to 1. Moreover,
now also $\sigma_{i/2}(e_j)$
will be a bounded net, converging strictly to 1.\\

Let $x$ be an element of $\mathscr{N}_\varphi\cap A$. Then
$xe_j\rightarrow x$ in norm. Moreover,
$\Lambda_\varphi(xe_j)=J\sigma_{i/2}(e_j)^*J\Lambda_\varphi(x)$.
Because $\sigma_{i/2}(e_j)$ also converges $^*$-strongly to 1, we
have $\Lambda_\varphi(xe_j)\rightarrow \Lambda_\varphi(x)$. Now if
$x$ is the norm-limit of $(a_i)$, with $a_i\in A_0$, then
$\Lambda_0(a_ie_j)=a_i\Lambda_\varphi(e_j)\rightarrow
x\Lambda_\varphi(e_j)=\Lambda_\varphi(xe_j)$ for each $e_j$. Since
$\Lambda_0$ is closed, each $xe_j$ and hence $x$ is in the domain of
$\Lambda_0$. So $\Lambda_0=\Lambda_\varphi$, and $A_0$ is a core for
$\Lambda_\varphi$.\\

The corollary follows, since the multiplicative unitary of $A$ and
the multiplicative unitary of $A_0$ on $\mathscr{H}_\varphi\otimes
\mathscr{H}_\varphi=\mathscr{H}_{\varphi_0}\otimes\mathscr{H}_{\varphi_0}$
coincide, and their first leg constitute respectively $A$ and the
C$^*$-algebraic quantum group associated to $A_0$.

\end{proof}

\noindent In our last proposition we will say something about the
dual of $(A_0,\Delta_0)$ when $A_0\subset A$ is a regular multiplier Hopf
algebra with $\varphi_{|A_0}\neq 0$.

\begin{Proposition}
Let $(\hat{A},\hat{\Delta})$ be the dual locally compact quantum
group of $(A,\Delta)$, and let $(\widehat{A_0},\widehat{\Delta_0})$
be the dual algebraic quantum group\footnote{The comultiplication in
$\widehat{A_0}$ is determined by
$\widehat{\Delta_0}(\omega)(x\otimes y)=\omega(yx)$. This is
\textit{not} the convention followed in \cite{VDae2}.} of
$(A_0,\Delta_0)$. Then
\[j: \widehat{A_0}\rightarrow \hat{A}: \varphi_0(\cdot\, a)
\rightarrow (\varphi(\cdot\, a)\otimes \iota)(W)\] is an injective
($^*$-)algebra homomorphism, such that \[(j\otimes
j)(\hat{\Delta}_0(\omega_1)(1\otimes\omega_2))=\hat{\Delta}(j(\omega_1))(1\otimes
j(\omega_2)),\] \[(j\otimes
j)(\hat{\Delta}_0(\omega_1)(\omega_2\otimes1))=\hat{\Delta}(j(\omega_1))(j(\omega_2)\otimes
1),\] for all $\omega_1,\omega_2\in \widehat{A_0}$.

\begin{proof}  Recall that $W$ denotes the multiplicative unitary of the left regular representation.
Remark that the expression
\[(\varphi(\cdot\, a)\otimes \iota)(W)\] makes sense, since
$\varphi(\cdot\, a)$ can also be written as a sum of elements of the
form $\varphi(b\cdot c)$, with $b,c$ in $A_0\subset
\mathscr{N}_\varphi\cap \mathscr{N}_\varphi^*$. It is also easily
seen that $j$ is injective.\\

\noindent We first check that $j$ preserves the $^*$-operation, in
case $A_0$ is a $^*$-algebraic quantum group. For simplicity, set
$\hat{a}=j(\varphi_0(\cdot\, a))$ for $a\in A_0$. Then
$\hat{a}\in\mathscr{D}(\Lambda_{\hat{\varphi}})$ and
$\Lambda_{\hat{\varphi}}(\hat{a})=\Lambda_{\varphi}(a)$, where
$\hat{\varphi}$ denotes the left Haar weight of the dual locally
compact quantum group $(\hat{A},\hat{\Delta})$. We know that
$\Lambda_\varphi(a)$ lies in the domain of the operator
$P^{1/2}J\delta^{-1/2}J$, where $P$ is the analytic generator for
the unitary one-parameter group
$\Lambda_\varphi(x)\rightarrow\Lambda_\varphi(\tau_t(x))$ and $J$ is
the modular conjugation for $\varphi$. But
$P^{1/2}J\delta^{-1/2}J=\hat{\nabla}^{1/2}$, where $\hat{\nabla}$ is
the modular operator for $\hat{\varphi}$. So $\hat{a}^*\in
\mathscr{D}(\Lambda_{\hat{\varphi}})$ and \begin{eqnarray*}
\Lambda_{\hat{\varphi}}(\hat{a}^*)&=&
\hat{J}\hat{\nabla}^{1/2}\Lambda_{\hat{\varphi}}(\hat{a})\\ &=&
\hat{J}P^{1/2}J\delta^{-1/2}J \Lambda_{\varphi}(a)\\ &=&
\hat{J}\Lambda_\varphi(\tau_{-i/2}(a)\delta^{-1/2})\\ &=&
\hat{J}\Lambda_\varphi((\tau_{-i/2}(a)\delta^{-1})\delta^{1/2})\\&=&
\hat{J}\Lambda_\psi(\tau_{-i/2}(a)\delta^{-1})\\
&=& \Lambda_\varphi(S(a)^*\delta)\\ &=&
\Lambda_\varphi(S_0(a)^*\delta_0) ,\end{eqnarray*} where we freely
used the formulas in \cite{Kus3}. Since $\varphi_0(\cdot
S_0(a)^*\delta_0)$ equals $\varphi_0(\cdot a)^*$, we arrive at $j(\varphi_0(\cdot a)^*)=j(\varphi_0(\cdot a))^*$.\\

\noindent Now we show that $j$ is an algebra morphism. Choose
$a,b\in A_0$. Choose $p_i,q_i$ in $A_0$ such that $a\otimes b=\sum
\Delta_0(p_i)(q_i\otimes 1)$. Then $\varphi_0(\cdot a)\cdot
\varphi_0(\cdot b)=\sum\varphi_0(q_i)\varphi_0(\cdot
p_i)=\sum\varphi(q_i)\varphi_0(\cdot p_i)$. It is enough then to
show that $\sum\varphi(q_i)\varphi(\cdot p_i)$ equals $\varphi(\cdot
a)\cdot \varphi(\cdot b)$ in $\mathscr{L}^1(A)$. But evaluating this
last functional in $x$, we get
$(\varphi\otimes\varphi)(\Delta(x)(a\otimes b))$, which equals
$\sum(\varphi\otimes \varphi)(\Delta(xp_i)(q_i\otimes
1))=\sum\varphi(q_i)\varphi(xp_i)$, so indeed both functionals are equal.\\

\noindent Similarly, $j$ respects the comultiplication. Namely:
\begin{eqnarray*}(\Lambda_{\hat{\varphi}}\otimes
\Lambda_{\hat{\varphi}})(\hat{\Delta}(\hat{a})(\hat{b}\otimes
1))&=&\Sigma W\Sigma
(\Lambda_{\hat{\varphi}}(\hat{a})\otimes\Lambda_{\hat{\varphi}}(\hat{b}))\\
&=& \Sigma W\Sigma (\Lambda_\varphi(a)\otimes\Lambda_\varphi(b))\\
&=& \Sigma
(\Lambda_\varphi\otimes\Lambda_\varphi)((S_0^{-1}\otimes\iota)(\Delta_0(b))
(a\otimes 1)), \end{eqnarray*} with $\Sigma$ denoting the flip. A
simple computation then shows that this expression is the image
under $(\Lambda_{\hat{\varphi}}\otimes \Lambda_{\hat{\varphi}})$ of
$(j\otimes j)(\hat{\Delta}_0(\varphi_0(\cdot\, a))(\varphi_0(\cdot\,
b)\otimes 1))$. Thus \[\hat{\Delta}(\hat{a})(\hat{b}\otimes
1)=(j\otimes
j)(\hat{\Delta}_0(\varphi_0(\cdot\,a))(\varphi_0(\cdot\,b)\otimes
1)).\] This implies the second equation of the proposition. The
first one follows similarly.

\end{proof}
\end{Proposition}

\noindent \textit{Remark.} The previous proposition says, that the
dual $\widehat{A_0}$ will be properly imbedded in $\hat{A}$ if $A_0$
is properly imbedded in $(A,\Delta)$. This implies that, under the
given conditions, also the dual $\widehat{A_0}$ of $A_0$ will have
an analytic structure.

\section{Structure of $^*$-algebraic quantum groups}

\noindent We apply the techniques of the above section to obtain
some interesting structural properties of $^*$-algebraic quantum
groups. While many of the results follow easily from the previous
section, we have decided to give new proofs, using only algebraic
machinery. As such, we can give a purely algebraic proof of the
existence of a \textit{positive} right invariant functional on a
$^*$-algebraic quantum group.\\

\noindent We fix a $^*$-algebraic quantum group $(A,\Delta)$ with
antipode $S$, positive left Haar functional $\varphi$, modular
automorphism $\sigma$ and modular element $\delta$. As a right Haar
functional (not assumed to be positive) we take $\psi=\varphi\circ
S$, with modular automorphism $\sigma'$. We adapt the proof of Lemma
\ref{lem3} to show that $A$ is spanned by eigenvectors for
$\kappa=\sigma^{-1}S^2$. We need a lemma.

\begin{Lemma} If $b$ is a non-zero element in $A_0$ and $n$ is an
even integer, then $b^*(\sigma'^nS^{2n})(b)\neq 0$.
\end{Lemma}
\begin{proof}

Suppose that $b\in A_0$ and $n\in 2\mathbb{Z}$ are such that
\[b^*(\sigma'^nS^{2n}(b))=0.\] Then
\[b^*\delta^n(\sigma^nS^{2n}(b))=0.\] Applying $\sigma^{-n/2}\circ S^{-n}$,
we get \[(\sigma^{n/2}S^n(b))^*\delta^n(\sigma^{n/2}S^n)(b)=0.\]
Since $\delta$ is self-adjoint, applying $\varphi$ to the previous
equation and using positivity and faithfulness of $\varphi$, we
obtain
\[\delta^{n/2}(\sigma^{n/2}S^n)(b)=0,\] hence $b=0$.

\end{proof}

\noindent \textit{Remark.} To be complete, we give an algebraic
argument demonstrating the self-adjointness of $\delta$. Choose $a$
and $b$ in $A_0$, then $\varphi(a^*a)b^*\delta b=(\varphi\otimes
\iota)((1\otimes b^*)\Delta(a^*a)(1\otimes
b))=\varphi(a^*a)b^*\delta^*b$. Applying $\varphi$ and using
polarization, we see that $\varphi(c^*\delta b)=\varphi(c^*\delta^*
b)$ for all $b,c\in A_0$. Finally, using the modular automorphism
and the fact that $A_0^2=A_0$, we get $\varphi(\delta
b)=\varphi(\delta^* b)$ for all $b\in A_0$. This implies $\delta
b=\delta^*b$ and $b\delta=b\delta^*$ for all $b\in A_0$. So
$\delta=\delta^*$.\\

\begin{Lemma}\label{lem5} If $a\in A$, then the linear span of the $\kappa^n(a)$, with $n\in \mathbb{Z}$,
is finite-dimensional.

\begin{proof}

We can follow the proof as in lemma \ref{lem3}:

Let $b$ be a fixed element of $A$.  Choose a non-zero $a\in A$, and
write \[a\otimes b = \sum_{i=1}^n \Delta(p_i)(1\otimes q_i), \] with
$p_i,q_i\in A$. Then
\[\kappa^n(a)\otimes \rho^{-n}(b)=\sum \Delta(p_i)(1\otimes \rho^{-n}(q_i)),
\qquad \textrm{for all }n\in \mathbb{Z},\] where $\rho=\sigma'S^2$.
Multiply this equation to the left with $1\otimes b^*$ to get
\[\kappa^n(a)\otimes b^*\rho^{-n}(b)=\sum
((1\otimes b^*)\Delta(p_i))(1\otimes \rho^{-n}(q_i)).\] Choose
$a_{ij},b_{ij}\in A$ such that \[(1\otimes
b^*)\Delta(p_i)=\sum_{j=1}^{m_i} a_{ij}\otimes b_{ij},\] and let $L$
be the finite-dimensional space spanned by the $a_{ij}$. We see that
$\kappa^n(a)\otimes b^*\rho^{-n}(b)\in L\odot A$, for every $n\in
\mathbb{Z}$. Using the previous lemma, we can conclude that
$\kappa^{2n}(a)\in L$ for all $n\in 2\mathbb{Z}$. But this easily
implies that the linear span $K$ of all $\kappa^n(a)$, with $n\in
\mathbb{Z}$, is a finite-dimensional, $\kappa$-invariant linear
subspace of $A$.

\end{proof}
\end{Lemma}

\noindent Denote by $(\hat{A},\hat{\Delta})$ the dual $^*$-algebraic
quantum group of $(A,\Delta)$. We can regard $\hat{A}$ and
$M(\hat{A})$ as functionals on $A$. We know from \cite{Kus5} that
$\hat{\delta}=\varepsilon \circ \kappa$ (this is also not so
difficult to prove algebraically). Then
\begin{eqnarray*} \langle \omega\hat{\delta},x\rangle &=& \langle \omega\otimes (\varepsilon\kappa),\Delta(x)\rangle
\\ &=& \langle \omega,\kappa(x)\rangle,\end{eqnarray*} for each
$\omega\in \hat{A}$ and $x\in A$. If $\omega$ is of the form
$\varphi(\cdot a)$, this means $\varphi(\cdot a)\hat{\delta}$ is a
scalar multiple of $\varphi(\cdot \kappa^{-1}(a))$. This implies
that for $\omega$ fixed, the linear span of the
$\omega\hat{\delta}^n$ is finite-dimensional. The same is of course
true for left
multiplication with $\hat{\delta}$.\\

\noindent By duality, we conclude that for each $a$ in $A$, the
linear span of the $\delta^na$ is a finite-dimensional space $K$.
(We could also prove this along the lines of Proposition 2.8.) Since
$\delta$ is a self-adjoint operator on $K$, with Hilbert space
structure induced by $\varphi$, we can diagonalize $\delta$. Hence
we arrive at

\begin{Proposition}\label{prop2} Let $(A,\Delta)$ be a $^*$-algebraic quantum group.
Then $A$ is spanned by elements which are eigenvectors for left
multiplication by $\delta$.
\end{Proposition}

\noindent We can use this to settle an open question (cf.
\cite{Kus1}):

\begin{Theorem}\label{the1} Let $(A,\Delta)$ be a $^*$-algebraic quantum
group. Then the scaling constant $\mu$ equals 1.
\begin{proof}

Choose a non-zero element $b\in A$ with $\delta b=\lambda b$, for
some $\lambda\in \mathbb{R}_0$. Then
$\varphi(bb^*\delta)=\lambda\varphi(bb^*)$. But the left hand side
equals $\mu\varphi(\delta bb^*) = \mu\lambda\varphi(bb^*)$. Since
$\varphi(bb^*)\neq 0$, we arrive at $\mu=1$.\\

\end{proof}

\end{Theorem}

\noindent Proposition \ref{prop2} can be strengthened:

\begin{Theorem}\label{the2} Let $(A,\Delta)$ be a $^*$-algebraic quantum group.
Then $A$ is spanned by elements which are simultaneously
eigenvectors for $S^2$, $\sigma$ and $\sigma'$, and left and right
multiplication by $\delta$. Moreover, the eigenvalues of these
actions are all positive.\end{Theorem}

\begin{proof} We know that $A$ is spanned by eigenvectors for left
multiplication with $\delta$, and the same is easily seen to be true
for $\kappa$ and $\rho=\sigma'S^2$. But all these actions commute.
Hence we can find a basis of $A$ consisting of simultaneous
eigenvectors. Since $\sigma, \sigma'$ and $S^2$ can be written as
compositions of the maps $\kappa,\rho$ and left and right
multiplication with $\delta$, the first part of the
theorem is proven.\\

We show that left multiplication with $\delta$ has positive
eigenvalues. This is easily done. Fix $a\in A_0$. If $\lambda$ is an
eigenvalue, choose an eigenvector $b$. Consider
$x=\Delta(a)(1\otimes b)$. Then $(\varphi\otimes \varphi)(x^*x)$
will be a positive number. But this is equal to
$\varphi(a^*a)\varphi(b^*\delta
b)=\lambda\varphi(a^*a)\varphi(b^*b)$. Hence $\lambda$ must be
positive. As before, duality implies that $\kappa$ and $\rho$ have
positive eigenvalues, hence the same is true of $\sigma, \sigma'$
and $S^2$.

\end{proof}

\noindent This theorem \textit{explains} why there exists an
analytic structure on a $^*$-algebraic quantum group $(A,\Delta)$:
the actions are all diagonal with positive entries! Hence $\sigma_z,
\sigma'_z, \tau_z$ and multiplication with $\delta^{iz}$
are all well-defined on $A$.\\

\noindent\label{right positivity}We can also see that
$\psi=\varphi\circ S$ is already a positive right invariant
functional, since
$\psi(a^*a)=\varphi(a^*a\delta)=\varphi((a\delta^{1/2})^*a\delta^{1/2})\geq
0$. Here we use that $\sigma(\delta^{1/2})=\delta^{1/2}$, which is easily proven using an eigenvector argument.\\

\noindent Finally remark that the extension of $\varphi$ to $M$,
with $M$ the von Neumann algebraic quantum group associated with
$A$, is an almost periodic weight, since the modular operator
$\nabla$ implementing $\sigma$ on $\mathscr{H}_\varphi$ is
diagonizable.

\section{Special cases}

\subsection*{Compact and discrete quantum groups}

\noindent Let $(A,\Delta)$ be a discrete locally compact quantum
group. Then $A$ is the C$^*$-algebraic direct sum of matrix algebras
$M_{n_\alpha}(\mathbb{C})$. The algebraic direct sum
$\mathscr{A}=\oplus_{\alpha}M_{n_\alpha}(\mathbb{C})$ has the
structure of a multiplier Hopf $^*$-algebra. So it is easy to see
that $\delta$, being a positive element in $\prod M_{n_\alpha}$, is
diagonizable with respect to $\mathscr{A}$. Then the same will be
true for $S^2$, the square of the antipode, since in a discrete
quantum group we have $S^2(a)=\delta^{-1/2}a\delta^{1/2}$. Lastly,
$\sigma$ is diagonizable since $\sigma=S^2$ in a discrete quantum
group.\\

\noindent Suppose now that $(A_0,\Delta_0)$ is a $^*$-algebraic
quantum group, properly imbedded in $(A,\Delta)$. Suppose $a$ is a
non-zero element in $A_0$ such that $a\notin \mathscr{A}$. We know
that $A_0$ has local units, so there exists $e\in A_0$ with $ae=a$.
Then $e^*a^*ae=a^*a$, and this implies that infinitely many
components of $e^*e$ have norm greater than 1. But this is
impossible,
since $e^*e\in A$. So $A_0\subset \mathscr{A}$.\\

\noindent The same argument implies that $A_0$ is again a
$^*$-algebraic quantum group of discrete type, since $A_0$ itself
will be an algebraic direct sum of matrix algebras. In particular,
$A_0$ has a co-integral $h_0$, which will be a grouplike projection
in $\mathscr{A}$. (A grouplike projection in a $^*$-algebraic
quantum group is a (self-adjoint) projection $p$ satisying
$\Delta(p)(1\otimes
p)=p\otimes p$. See \cite{Lan2} for more details.)\\

\noindent  The dual side is also easy to treat. Namely, let
$(A,\Delta)$ be a (reduced) compact locally compact quantum group.
We know then that $A$ contains a dense multiplier Hopf $^*$-algebra
$\mathscr{A}$. Suppose that $(A_0,\Delta_0)$ is a multiplier Hopf
$^*$-algebra imbedded in $(A,\Delta)$. Since the left Haar weight
$\varphi$ is everywhere defined, the elements of $A_0$ are
automatically integrable. Then $(A_0,\Delta_0)$ will be a
$^*$-algebraic quantum group. We know that $\widehat{A_0}$ is a
discrete quantum group properly imbedded in $\hat{A}$. Hence
$A_0\subset \mathscr{A}$ and $(A_0,\Delta_0)$ is a compact
$^*$-algebraic quantum group. The dual $p$ of the co-integral $h_0$
of $\widehat{A_0}$ in $\hat{\mathscr{A}}$ will be a grouplike
projection in $\mathscr{A}$. It will be a unit for $A_0$.

\subsection*{Locally compact groups}

\noindent Suppose $G$ is a locally compact group. Let
$(A_0,\Delta_0)$ be a regular multiplier Hopf $^*$-algebra imbedded
in $(\mathscr{L}^\infty(G),\Delta)$, where $\Delta$ is the usual
comultiplication determined by $\Delta(f)(g,h)=f(gh)$. Then
$A_0\subset M(C_0(G))=C_b(G)$, so $A_0$ consists of bounded
continuous functions on $G$. Let $\bar{A_0}$ be the normclosure of
$A_0$ in $\mathscr{L}^\infty(G)$. Then $\Delta$ restricts to a
$^*$-algebra morphism $\bar{A_0}\rightarrow M(\bar{A_0}\otimes
\bar{A_0})$. Since $\bar{A_0}$ is abelian, this induces a locally
compact semigroup structure on the spectrum $X$ of $\bar{A_0}$.
Since $S_0$ extends to an $S$ on $\bar{A_0}$, the semigroup will be
a locally compact group. But this means that $X$ has a Haar measure.
So $(A_0,\Delta_0)$ is properly imbedded in the C$^*$-algebraic
quantum group $(C_0(X),\Delta)$, hence is itself a $^*$-algebraic
quantum group.\\

\noindent Remark however that the Haar functional on $A_0$ can be
different from integration with respect to the Haar measure on $G$.
Consider for example the linear span $A_0$ of the functions
$f_x:t\rightarrow e^{itx}$ in $\mathscr{L}^\infty(\mathbb{R})$, with
$x\in \mathbb{R}$. Then $A_0$ is a compact Hopf $^*$-algebra, but
none of its non-zero elements are integrable with respect to the
Lebesgue measure. It is easy to see that the Haar functional
$\varphi_0$ on $A_0$ is given by $\varphi_0(f_x)=\delta_{x,0}$, and
that the space $X$ equals $\mathbb{R}$ with the discrete
topology.\\

\noindent The dual case is not so clear: suppose $G$ is a locally
compact group, and $(A_0,\Delta_0)$ is imbedded in
$(\mathscr{L}(G),\Delta)$, where $\Delta$ is determined by
$\Delta(u_g)=u_g\otimes u_g$ on the generators of $\mathscr{L}(G)$.
Will the $C^*$-algebraic closure $\bar{A_0}$ with the restriction of
$\Delta$ be of the form $(C^*(X),\Delta)$ for some locally compact
group $X$? This will of course be true if $A_0$ is properly imbedded
in $C^*_r(G)$, since then we can apply the theory of the second
section to conclude that $\bar{A_0}$ is a cocommutative
C$^*$-algebraic quantum group with normdense $^*$-algebraic quantum
group, hence of the form $(C^*_r(X),\Delta)$.\\

\noindent Let us now look at the results of the third section in the
commutative case. Let $G$ be a locally compact group with a compact
open subgroup $H$. Consider the regular functions on $H$ - i.e. the
functions generated by the matrix-coefficients of finite-dimensional
irreducible representations of $H$. We can see them as functions on
$G$. The linear span of left translates of these functions by
elements of $G$ is denoted by $P_0(G)$. In \cite{Lan1}, it is shown
that $P_0(G)$ forms a dense multiplier Hopf$^*$-algebra inside
$(C_0(G),\Delta)$, with the usual comultiplication, and that every
commutative $^*$-algebraic quantum group is of this form. In this
setting, the only non-trivial object is the modular function
$\delta$. According to our results, it should be diagonizable. This
is easily seen to be true. For example, the characteristic function
of $H$ will be an eigenvector for left multiplication. Indeed: the
Haar measure on $H$ is the restriction of the Haar measure on $G$.
Hence $\delta_{|H}$ is the modular function of $H$. Since $H$ is
compact, $\delta_{|H}=1$. So every regular function on $H$ is
invariant for left multiplication. Then the translates by some
element $g$ of such functions will be eigenvectors with eigenvalue
$\delta(g)$, and the
linear span of all such translates equals $P_0(G)$.\\

\subsection*{The case of the quantum groups $U_q(su(2))$ and
$SU_q(2)$}

\noindent Finally, we consider a particular, non-trivial example of
a multiplier Hopf $^*$-algebra $(A_0,\Delta_0)$, imbedded in the
multiplier algebra of a discrete $^*$-algebraic quantum group
$(\mathscr{A},\Delta)$. This is not a situation we have discussed,
since this multiplier algebra contains unbounded operators (when
acting on the Hilbert space closure of $\mathscr{A}$ by left
multiplication). We will see which of our results are still true in
this case.\\

\noindent So as $(A_0,\Delta_0)$, we take the quantum enveloping Lie
algebra $U_q(su(2))$, with $q$ nonzero in $\rbrack -1,1\lbrack$. As
a $^*$-algebra, it is generated by two elements, $E$ and $K$, with
$K$ invertible and self-adjoint, obeying the following commutation
relations:
\begin{displaymath} \left\{ \begin{array}{l} EK = q^{-1}KE\\
 \lbrack E, E^* \rbrack = \frac{1}{q-q^{-1}}(K^2-K^{-2}). \end{array}
\right. \end{displaymath} The comultiplication on the generators is
given by
\begin{displaymath} \left\{
\begin{array}{l} \Delta_0(K) = K \otimes K \\ \Delta_0(E) = E \otimes K + K^{-1} \otimes E.\end{array} \right.
\end{displaymath} To see that this comultiplication is
well-defined, it is enough to check that it respects the commutation
relations, but this is easily done. The antipode is determined by
\begin{displaymath} \left\{
\begin{array}{l} S_0(K)= K^{-1} \\ S_0(E) = -qE \\
S_0(E^*)= -q^{-1}E^*.\end{array} \right.
\end{displaymath}\\

\noindent As our $^*$-algebraic quantum group
$(\mathscr{A},\Delta)$, we take the $^*$-algebraic quantum group
$\widehat{\mathscr{B}}$, where $\mathscr{B}$ is the compact
$^*$-algebraic quantum group associated with  $SU_q(2)$, Woronowicz'
twisted $SU(2)$-group. As a $^*$-algebra, $\mathscr{B}$ is generated
by two elements $a$ and $b$, such that
\[ \left\{ \begin{array}{ll}
ab=qba \\ ab^*=qb^*a\\
\lbrack b,b^*\rbrack =0 \\
a^*a = 1-q^{-2}b^*b \\
aa^* = 1-b^*b. \\
\end{array} \right. \]The co-multiplication is given by: \[\left\{ \begin{array}{ll}
\Delta (a) = a \otimes a -q^{-1} b \otimes b^* \\
\Delta (b) = a \otimes b + b \otimes a^*. \end{array}\right.\] For
convenience, we state that the antipode $S$ is given by
\[\left\{
\begin{array}{l} S(a) = a^*\\ S(a^*)=a\\S(b)=-q^{-1}b\\ S(b^*)=-qb^*.\end{array}\right.\]
We will not need the concrete description of the left invariant
functional, but we need to know the modular group, which we now
denote by $\rho$. To be complete, we also provide the scaling group,
which we will denote by $\theta$:
\[\left.\begin{array}{ll}\left\{\begin{array}{l}
\rho_z(a)=q^{-2iz}a\\ \rho_z(b)=b\end{array}\right. &
\left\{\begin{array}{l} \theta_z(a)=a\\
\theta_z(b)=q^{-2iz}b.\end{array}\right.\end{array}\right.\] The
modular element will of course be trivial, since the quantum group
is compact.

\noindent The easiest way to see that $A_0$ can be imbedded in
$M(\widehat{\mathscr{B}})$, is by creating a pairing\footnote{In
this part we follow the algebraic convention $\langle
\Delta(b),x\otimes y\rangle = \langle b,xy\rangle$.} between
$\mathscr{B}$ and $A_0$. For, since $\mathscr{B}$ is compact, it is
known that $M(\widehat{\mathscr{B}})$ can be identified with the
vector space of \textit{all} linear functionals on $\mathscr{B}$.
The fact that there is a \emph{pairing}, implies that the inclusion
of $A_0$ in $M(\widehat{\mathscr{B}})$ will be a morphism. The
concrete pairing is as follows:
\[\left.\begin{array}{ll}\left\{ \begin{array}{l}  \langle K,a \rangle = q^{-1/2}\\ \langle K,a^* \rangle =q^{1/2}\\
\langle K,b \rangle = 0 \\ \langle K,b^* \rangle  =
0\end{array}\right. & \left\{\begin{array}{l}
\langle E,a \rangle = 0 \\ \langle E,a^* \rangle = 0. \\ \langle E,b \rangle  = 0 \\  \langle E,b^* \rangle = -q.  \\
\end{array}\right.\end{array}\right.\]\\

\noindent Since on the dual of a compact algebraic quantum group the
modular group $\sigma$ and the scaling group $\tau$ coincide, we
find the following behavior of $A_0$:
\[\left\{\begin{array}{ll} \sigma_z(K)=\tau_z(K)=K\\
\sigma_z(E)=\tau_z(E)=q^{2iz}E.\end{array}\right.\] But although
there is general invariance under the scaling (and thus the modular)
group, we no longer have that $A_0$ is invariant under left
multiplication by $\delta^{z}$, with $\delta$ the modular element of
$\widehat{\mathscr{B}}$. For this would imply that actually
$\delta^{z}\in A_0$, since $A_0$ is a Hopf algebra. Remark that this
$\delta^{z}$ is easily computable, for it is given as a functional
by $\varepsilon\circ \rho_{iz}$, with $\varepsilon$ the co-unit of
$\mathscr{B}$. We find that applying $\delta$ is the same as pairing
with $K^{-4}=(K^*K)^{-2}$, so uniqueness gives us that
$\delta^{it}=K^{-4it}$. It is clear that this is no element in
$A_0$. Remark also, that right or left multiplication with $\delta$
is no longer diagonal. This is easy to see, using that $A_0$ has
$\{K^lE^mF^n\mid l\in \mathbb{Z},m,n\in \mathbb{N}\}$ as a basis. In
fact, since span$\{K^{4n}X\}$ has infinite dimension for any $X\in
A_0$, we get that $A_0\cap \widehat{\mathscr{B}}=\{0\}$.
\\

\noindent We note that in this example we are in a special
situation: $(\widehat{SU_q(2)},\hat{\Delta})$ is the C$^*$-algebraic
quantum group generated by $K$, $K^{-1}$ and $E$, in the sense of
Woronowicz. Moreover, the multiplier Hopf $^*$-subalgebra is linked
by a pairing to a $^*$-algebraic quantum group. This could explain
why we still have invariance under $\tau_t$ and $\sigma_t$. For
example, the same type of behavior occurs with the quantum
$az+b$-group. Remark that in these cases, the corresponding Hopf
$^*$-algebra can be viewed as the infinitesimal version of the
quantum group. We do not know if it is a general fact that the
one-parameter groups descend to the Hopf $^*$-algebra associated
with the quantum group, if such an object is present. In any case,
the connection between a locally compact quantum group and a Hopf
$^*$-algebra representing the quantum group at an infinitesimal
level, is at present not well understood in a general framework.
\\

\end{document}